\documentclass[a4paper, 11pt]{amsart}
\usepackage{amsmath, amssymb, amsthm, amsfonts,enumitem,tikz,relsize}
\usepackage[normalem]{ulem}
\usepackage[colorlinks = true, citecolor = blue, hypertexnames=false]{hyperref}
\usepackage{xcolor}
\hypersetup{
    colorlinks,
    linkcolor={red!70!black},
    citecolor={blue!70!black},
    urlcolor={blue!80!black}
}
\usepackage{framed,tabu,float,caption,scalerel}
\usepackage[OT1]{fontenc}
\usepackage{comment}

\theoremstyle{plain}
\newtheorem*{maintheorem}{\hypertarget{h:theorem}{Theorem}}
\newtheorem*{maincorollary}{\hypertarget{h:corollary}{Corollary}}
\newtheorem{theorem}{Theorem}[section]
\newtheorem{corollary}[theorem]{Corollary}
\newtheorem{lemma}[theorem]{Lemma}

\newtheorem{fact}[theorem]{Fact}
\newtheorem*{corollary*}{Corollary}

\theoremstyle{definition}
\newtheorem*{definition*}{Definition}
\newtheorem*{example*}{Example}

\newtheorem*{problem*}{Problem}

\newtheorem*{conjecture*}{Conjecture}
\newtheorem*{mainconjecture}{\hypertarget{h:conjecture}{Conjecture}}

\theoremstyle{remark}
\newtheorem{claim}{Claim}

\newenvironment{proofclaim}[1][Proof of Claim]{\begin{proof}[#1]}{\end{proof}}

\DeclareMathOperator{\rk}{rk}
\DeclareMathOperator{\gtd}{gtd}

\DeclareMathOperator{\Aut}{Aut}

\DeclareMathOperator{\End}{End} 
\DeclareMathOperator{\Sym}{Sym} 
\DeclareMathOperator{\Alt}{Alt}

\DeclareMathOperator{\GL}{GL}
 
\DeclareMathOperator{\AGL}{AGL} 
\DeclareMathOperator{\PGL}{PGL}

\newcommand{\bK}{\mathbb{K}}

\begin{document}
\title[Actions of $\Alt(n)$ on groups of finite Morley rank]{Actions of $\Alt(n)$ on groups of finite Morley rank without involutions}
\author{Tuna Alt\i{}nel}
\address{Universit\'e de Lyon \\
Universit\'e Claude Bernard Lyon 1\\
CNRS UMR 5208\\
Institut Camille Jordan\\
43 blvd du 11 novembre 1918\\
F-69622 Villeurbanne cedex\\
France}
\email{altinel@math.univ-lyon1.fr}
\author{Joshua Wiscons}
\address{Department of Mathematics and Statistics\\
California State University, Sacramento\\
Sacramento, CA 95819, USA}
\email{joshua.wiscons@csus.edu}
\thanks{The work of the second author was partially supported by the National Science Foundation under grant No. DMS-1954127.}
\date{\today}

\begin{abstract}
We investigate faithful representations of $\Alt(n)$ as automorphisms of a connected group $G$ of finite Morley rank. We target a lower bound of $n$ on the rank of such a \emph{nonsolvable} $G$, and our main result achieves this in the case when $G$ is  without involutions. In the course of our analysis, we also prove a corresponding bound for solvable $G$ by leveraging recent results on the abelian case. We conclude with an  application towards establishing natural limits to the degree of generic transitivity for permutation groups of finite Morley rank.  
\end{abstract}

\maketitle

% % % % % % % % % % % % % % % % % % % % % % % % % % % % % % % % % % % % % % % % % % % % % % % % % % % % % % % % %
% % % % % % % % % % % % % % % % % % % % % % % % % % % % % % % % % % % % % % % % % % % % % % % % % % % % % % % % %
% % % % % % % % % % % % % % % % % % % % % % % % % % % % % % % % % % % % % % % % % % % % % % % % % % % % % % % % %
% Introduction 
% % % % % % % % % % % % % % % % % % % % % % % % % % % % % % % % % % % % % % % % % % % % % % % % % % % % % % % % %
% % % % % % % % % % % % % % % % % % % % % % % % % % % % % % % % % % % % % % % % % % % % % % % % % % % % % % % % %
% % % % % % % % % % % % % % % % % % % % % % % % % % % % % % % % % % % % % % % % % % % % % % % % % % % % % % % % %
\section{Introduction}\label{S.introduction}

Morley rank is a fundamental, model-theoretic notion of dimension generalizing Zariski dimension from algebraic geometry. The affine algebraic groups over algebraically closed fields are among the groups of finite Morley rank as are various nonalgebraic groups. The importance of groups of finite Morley rank is highlighted by---and rooted in---their appearance as ``binding groups'' (analogous to Galois groups) in the analysis of certain model-theoretic classes of structures, and this connection makes the study of \emph{permutation groups} of finite Morley rank  particularly salient. 
The recent work  of James Freitag and Rahim Moosa---later refined in collaboration with R\'{e}mi Jaoui---illustrates this point  well where  bounding the degree of generic transitivity (as defined in Section~\ref{S.PermGroups}) of certain binding groups yields new information about algebraic relationships between nonalgebraic solutions to %algebraic
differential equations \cite{FrMo22,FJM22}.

Much of the current work on permutation groups of finite Morley rank is structured around  guiding problems of Alexandre Borovik and Gregory Cherlin about the degree of {generic} transitivity~\cite{BoCh08}. This topic originated in the algebraic setting in work of Vladimir Popov~\cite{PoV07}, and his findings---as well as   the questions that remain---have provided much inspiration for the finite Morley rank setting.  And in the other direction, results obtained in this more general (and perhaps more appropriate) context have also extended the algebraic theory. 

Due in large part to work on generic transitivity, there has been much  recent activity investigating \emph{representations} in a finite Morley rank context, and one emerging theme is the relevance of understanding representations  of \emph{finite groups} in this setting. This was already present in the   work of Borovik and Cherlin where the determination of the minimal Morley rank of an abelian group on which $\Sym(n)$ can act faithfully and definably drives their subsequent analysis establishing a bound on the degree of generic transitivity for primitive groups of finite Morley rank. The authors of this article also studied representations of the symmetric group in their classification of generically $4$-transitive actions on sets of rank $2$~\cite{AlWi19}, and in that work, there was a need to analyze representations where the acted upon group was not necessarily abelian.  Additionally, the recent work of Ay\c{s}e Berkman and Alexandre Borovik completing the classification of representations of finite Morley rank with a sufficiently high degree of generic transitivity makes critical use of the hyperoctahedral group~\cite{BeBo18,BeBo21}. 

The present authors have sought to extend their analysis of actions on sets of  rank $2$ to sets of arbitrary rank, and additional knowledge of representations of the symmetric group appears to be essential. 
In particular, a lower bound on the rank of a (not necessarily abelian) group admitting a faithful and definable action of $\Sym(n)$ is strongly desired, and a classification of those actions achieving the lower bound would further clarify the situation. The latter was recently addressed by Luis Jaime Corredor, Adrien Deloro, and the second author~\cite{CDW21} (with the work of Borovik on actions of arbitrary finite groups being highly relevant as well~\cite{BoA20}); the following conjecture proposes a rough bound to address the former.

\begin{mainconjecture}[see {\cite[Section~1.2]{CDW21}}]
Suppose that $\Sym(n)$ or $\Alt(n)$ acts definably and faithfully by automorphisms on a connected group of finite Morley rank $G$, with $n$ sufficiently large. If $G$ is nonabelian, then $\rk G \ge n-1$, and if $G$ is nonsolvable, then $\rk G \ge n$.
\end{mainconjecture}

The solvable case follows readily from the main result of~\cite{CDW21}; we address this first (see Lemma~\ref{lem.solvable}). The nonsolvable case is much less clear to us.   When $G$ is algebraic, routine reduction to $G$ being simple leads to an embedding of $\Alt(n)$ into $G$, so this appears to be under control, if not already classical. However, small values of $n$ are indeed complicated by exceptional isomorphisms, for both $\Sym(n)$ and $\Alt(n)$. In this paper, we establish the \hyperlink{h:conjecture}{Conjecture} at the other (highly nonalgebraic) extreme. 

\begin{maintheorem}
Suppose  $\Alt(n)$, with $n\ge 4$, acts definably and faithfully on a connected nonsolvable group $G$ of finite Morley rank without involutions. Then $\rk G \ge n$.
\end{maintheorem}

In addition to making progress on the \hyperlink{h:conjecture}{Conjecture}, the \hyperlink{h:theorem}{Theorem} supplements the rather small amount of existing work controlling the structure of nonsolvable groups of finite Morley rank without involutions. According to the Algebraicity Conjecture of Gregory Cherlin and Boris Zilber, no such group should exist, but this appears to be quite far from being resolved. Indeed, it was not until 2018 that Oliver Fr\'econ established the Algebraicity Conjecture for groups of rank $3$ \cite{FrO18}, and most other work on nonsolvable groups without involutions is confined to special configurations such as  ``bad groups.'' 

As mentioned before, the \hyperlink{h:theorem}{Theorem}  has applications to permutation groups with a high degree of generic transitivity. This is discussed in some detail in Section~\ref{S.PermGroups} where we also explicitly state our \hyperlink{h:corollary}{Corollary} for such actions.

% % % % % % % % % % % % % % % % % % % % % % % % % % % % % % % % % % % % % % % % % % % % % % % % % % % % % % % % %
% % % % % % % % % % % % % % % % % % % % % % % % % % % % % % % % % % % % % % % % % % % % % % % % % % % % % % % % %
% SECTION 
% % % % % % % % % % % % % % % % % % % % % % % % % % % % % % % % % % % % % % % % % % % % % % % % % % % % % % % % %
% % % % % % % % % % % % % % % % % % % % % % % % % % % % % % % % % % % % % % % % % % % % % % % % % % % % % % % % %
\section{The Theorem}

Our proof of the \hyperlink{h:theorem}{Theorem} is inductive. We first collect a handful of background results. The first two (Facts~\ref{fact.SmallGroups} and \ref{fact.InvolutiveAut}) allow us to quickly dispense with the tightest, low rank configurations,  and the next two (Facts~\ref{fact.4Group} and \ref{fact.BigCorank}) drive our inductive approach. Ultimately, we are forced to consider the case when the acted upon group is solvable; we treat this quickly in Lemma~\ref{lem.solvable} by leveraging recent work on the abelian case (where the context is significantly more general than finite Morley rank and the conclusion significantly stronger than just a bound on the rank)~\cite{CDW21}.

General background on groups of finite Morley rank can be found in \cite{ABC08} or \cite{BoNe94}, among  other sources. For the most part, the reader familiar with the algebraic setting but not the logical one may safely replace \emph{definable}  with {constructible} and \emph{rank} with {dimension}.

% % % % % % % % % % % % % % % % % % % % % % % % % % % % % % % % % % % % % % % % % % % % % % % % % % % % % % % % %
% Subsubsection 
% % % % % % % % % % % % % % % % % % % % % % % % % % % % % % % % % % % % % % % % % % % % % % % % % % % % % % % % %
\subsection{Background and preliminary results}
Our first two background results use the term \emph{bad group}, which we here take to mean a nonsolvable connected group of finite Morley rank all of whose proper definable connected subgroups are nilpotent. 

\begin{fact}[\cite{ReJ75,ChG79,WiJ16,DeWi16,FrO18}]\label{fact.SmallGroups}
If $G$ is a connected nonsolvable group of finite Morley rank without involutions, then $\rk G \ge 4$. Moreover, if $\rk G\le 5$, then $G'$ is quasisimple and bad.
\end{fact}

In the presence of involutive automorphisms, the conclusion of Fact~\ref{fact.SmallGroups} is further clarified by the following classical fact.

\begin{fact}[Delahan-Nesin, see~{\cite[Proposition~13.4]{BoNe94}}]\label{fact.InvolutiveAut}
A simple bad group of finite Morley rank has no definable involutive automorphism.
\end{fact}

The next two facts  underlie our inductive approach to the \hyperlink{h:theorem}{Theorem}.

\begin{fact}[{\cite[Theorem~5]{BBC07}}]\label{fact.4Group}
Let $G$ be a connected group of finite Morley rank without involutions and $K$ an elementary abelian group of order $4$ acting definably on $G$. Then $G = \langle C_G^\circ(\alpha) \mid \alpha\in K\setminus\{1\}\rangle$.
\end{fact}

We mention in passing a related result about actions of $8$-groups; it was a part of our original proof of the \hyperlink{h:theorem}{Theorem} but has since been removed (at the cost of introducing several more lines of mathematics). The result---which is due to Jeffrey Burdges, though not in any of his papers---is that if an elementary abelian group $E$ of order $8$ acts definably on a group of finite Morley rank $G$ without involutions, then $C_G^\circ(e)$ is nonsolvable for some involution $e\in E$.

The next result is an adaptation of \cite[Proposition~2.3]{BoCh08} (together with \cite[Lemma~2.4]{BoCh08}). %; see also  \cite[Lemma~4.23]{AlWi18}. 
We only need the corollary that follows for our proof of the \hyperlink{h:theorem}{Theorem}, so 
we delay the definition of generically $t$-transitive actions until Section~\ref{S.PermGroups}.

\begin{fact}[see {\cite[Proposition~2.3]{BoCh08}}]\label{fact.BigCorank}
Let $G$ be a connected group of finite Morley rank and $H$ a maximal definable connected subgroup. Assume  $H$ has a finite core, i.e.~$\bigcap_{g \in G} H^g$ is finite. Let $X$ be the set of right cosets of $H$, and let  $\gtd(G,X)$ denote the maximum $t$ such that the action of $G$ on $X$ is generically $t$-transitive. Setting $s = \rk X$, we have  
\[\rk G \le s\cdot \gtd(G,X) + s(s-1)/2.\]
\end{fact}

\begin{corollary}\label{cor.BigCorank}
Let $G$ be an infinite simple group of finite Morley rank without involutions and $H$ a proper definable connected subgroup of corank~$s$. Then $\rk G \le s(s+1)/2$ and 
$\rk H \le \rk G+\frac{1}{2}-\sqrt{2\rk G+\frac{1}{4}}.$
\end{corollary}
\begin{proof}
We may assume $H$ is a maximal definable connected subgroup. Since $G$ is without involutions, no action of $G$ on an infinite set is generically $2$-transitive; indeed, generic $2$-transitivity would imply $G$ contains an element swapping a generic pair from the set being acted upon, forcing $G$ to contain an involution.   Fact~\ref{fact.BigCorank} (with $\gtd(G,X)=1$) now yields the desired result.
\end{proof}

Asymptotically, Corollary~\ref{cor.BigCorank} says little, but it does create a little space for induction and will be used in the sequel.

Finally, we address the solvable analogue of the \hyperlink{h:theorem}{Theorem}. As expected, it follows quickly from the recent solution of the abelian case, which we state here in a very abbreviated form. When $n\ge 7$, our formulation follows directly from the Theorem of \cite{CDW21}; when $n< 7$, it is a combination of the 
First Geometrisation Lemma and the Recognition Lemma of \cite{CDW21}.

\begin{fact}[see {\cite{CDW21}}]\label{fact.CDW}
Let $n\ge 3$. If $V$ is a faithful connected $\Sym(n)$-module of finite Morley rank, then $\dim V \ge n-2$, and if equality holds, then $V$ is an elementary abelian $p$-group for some prime $p$ dividing $n$. 

Moreover, the same is true for $\Alt(n)$-modules provided either $n \geq 10$ or both  $n\ge 7$ and $V$ is $2$-divisible.
\end{fact}

Our statement of the solvable case requires a definition.

\begin{definition*}
Let $G$ be a connected group of finite Morley rank, and suppose a group $H$ acts definably on $G$. We say that $G = G_\ell > G_{\ell-1} >\cdots > G_0 = 1$ is an \emph{$H$-minimal} series for $G$ if for all $i>0$,
\begin{itemize}
    \item $G_i$ is definable, nontrivial, connected, normal, and $H$-invariant; and
    \item each $G_i/G_{i-1}$ possesses \emph{no} proper definable, nontrivial, connected, normal, $H$-invariant subgroups.
\end{itemize}
\end{definition*}

\begin{lemma}[Solvable Lemma]\label{lem.solvable}
Let $n\ge 3$. Suppose $\Sym(n)$ acts faithfully and definably on a connected solvable group $B$ of finite Morley rank. If $B = B_\ell > B_{\ell-1} >\cdots > B_0 = 1$ is any $\Sym(n)$-minimal abelian series, then $\rk B \ge n + \ell - 3$; in particular, $\rk B \ge n -1$ whenever $B$ is nonabelian.

Moreover, the same is true for an action of $\Alt(n)$ provided either $n \geq 10$ or both  $n\ge 7$ and $B$ is $2$-divisible.
\end{lemma}
\begin{proof}
Let $S := \Sym(n)$, $A := S'$, and $r = \rk B$. There are $\ell$ quotients in the series, each of rank at least $1$, so $s_m := \rk\left (B_m/B_{m-1}\right)\le  r - \ell + 1$. 

Towards a contradiction, assume $r < n + \ell - 3$ (which also implies $n\ge 4$). Then for all $1\leq m\leq l$, we have $s_m \le r - \ell + 1 < n-2$, and Fact~\ref{fact.CDW} implies that $A$ does \emph{not} act faithfully on any $B_m/B_{m-1}$. 

If $n\ge 5$, then $A$ is simple and $[A,\ldots,A,B] = 1$, and (as $A$ is simple) this easily implies $[A,B] = 1$, contradicting faithfulness. It remains to consider when the acting group is $S = \Sym(4)$ in which case $s_m = 1$ for all $m$. It is easily seen that each involution in $S$ must either centralize or invert any connected rank $1$ section on which it acts (see  \cite[I,~Lemma~10.3]{ABC08}), so in this case we again find that $[A,\ldots,A,B] = 1$. This condition, for a faithful action, in fact implies that $A$ is nilpotent (see for example \cite{HaP58}), again a contradiction.
\end{proof}

% % % % % % % % % % % % % % % % % % % % % % % % % % % % % % % % % % % % % % % % % % % % % % % % % % % % % % % % %
% Subsubsection 
% % % % % % % % % % % % % % % % % % % % % % % % % % % % % % % % % % % % % % % % % % % % % % % % % % % % % % % % %
\subsection{The proof}
Our main analysis starts here. We first deal with the cases of lowest rank.

\begin{lemma}\label{lem.KleinOnSmall}
Suppose a Klein four-group $K$ acts definably and faithfully on a connected nonsolvable group $G$ without involutions. Then $\rk G \ge 6$. 
\end{lemma}
\begin{proof}
In light of Fact~\ref{fact.SmallGroups}, we need only consider when $4\le\rk G\le 5$ and $G'$ is quasisimple and bad, which further implies $\rk(G')\ge 4$. Crucially, Fact~\ref{fact.InvolutiveAut} tells us that $[K, G'] \le Z(G')$, and as $G'$ is connected, we in fact have that $[K, G'] \le Z^\circ(G')$ (see \cite[Corollary~5.29]{BoNe94}).

First assume $G' = G$. We study the action of $K$ on the series $1\le Z^\circ(G) < G$, and we already know $[K, G] \le Z^\circ(G)$. By Fact~\ref{fact.SmallGroups}, $\rk(G/Z(G)) \ge 4$, so $\rk Z^\circ(G) \le 1$. Consequently, each involution in $K$  either centralizes or inverts $Z^\circ(G)$, so some involution $\alpha \in K$ must centralize $Z^\circ(G)$. Thus, $\alpha$ acts quadratically on $G$ (i.e.~ $[\alpha,\alpha,G] = 1$), which implies  $[\alpha,G]$ has exponent $2$. As $G$ is without involutions, $[\alpha,G]=1$, contradicting faithfulness. 

Next suppose that $G' \neq G$. We now consider  $1< G' < G$. Then $G'$ has rank $4$, and Fact~\ref{fact.SmallGroups} forces $Z(G')$ to be finite. So here we have $[K, G'] = 1$, and consideration of the rank $1$ quotient $G/G'$ shows (as before) that some $\alpha \in K$ acts quadratically on $G$, a contradiction (as before).
\end{proof}

\begin{proof}[Proof of the \hyperlink{h:theorem}{Theorem}]
Assume the theorem is false. Choose $n$ minimal such that the theorem fails, and then choose $G$ witnessing this failure that minimizes $\rk G$ among all such witnesses. 

Set $A :=  \Alt(n)$ and $r:=\rk G$. As $n\ge 4$, $A$ contains a Klein $4$-group, so $r\ge 6$ by Lemma~\ref{lem.KleinOnSmall}. As we are working with a counterexample, $n\ge 7$.

\begin{claim}
We may assume $G$ is simple.
\end{claim}
\begin{proofclaim}
Let  $G = G_\ell > G_{\ell-1} >\cdots > G_0 = 1$ be an $A$-minimal series for $G$; thus, each factor $G_i/G_{i-1}$ is either abelian or nonsolvable. 

We show $\ell = 1$. Assume not; then $\rk(G_i/G_{i-1}) < r$ for each $i$. Moreover, if a factor $G_i/G_{i-1}$ is abelian, then $\rk(G_i/G_{i-1}) < r-2$ as otherwise $G$ would be solvable (since connected groups of rank at most $2$ are solvable). Consequently, Fact~\ref{fact.CDW} ensures that each abelian factor of the series (which is $2$-divisible by our assumption on $G$) is centralized by $A$, and minimality of $r$ ensures that each nonsolvable factor is also centralized by $A$. Thus, in this case, $[A,\ldots,A,G] = 1$, so as $A$ is simple, $[A,G] = 1$, a contradiction.  

We now show $G$ is quasisimple. As $\ell =1$, $G$ is equal to its generalized Fitting subgroup (which also equals the socle of $G$), so there exist definable, normal, connected, quasisimple subgroups $Q_1,\ldots,Q_k$ of $G$ such that $G=Q_1\oplus\ldots\oplus Q_k$.
All $Q_i$ are connected, so $k\le r < n$ (in fact, $4k\le r<n$ by Fact~\ref{fact.SmallGroups}). Since $A$ acts trivially on sets of size less than $n$, $A$ normalizes each $Q_i$, so as $\ell =1$, $G = Q_1$ is quasisimple.

Again using that $\ell =1$,  $G$ is quasisimple with \emph{finite} center. Since $[A,G]$ is connected, $A$ must act faithfully on $G/Z(G)$ as otherwise $A$ would centralize all of $G$. Hence, we may replace $G$ with this quotient.
\end{proofclaim}

For $\sigma,\psi\in A$, we use the notation $C_\sigma := C_G^\circ(\sigma)$ and $C_{\sigma,\psi} := C_G^\circ(\sigma,\psi)$. We  write $A_{\sigma}$ to denote the subgroup of permutations whose support is \emph{contained} in the support of $\sigma$;  $A_{\sigma^\perp}$ denotes the subgroup of permutations with support \emph{disjoint} from the support of $\sigma$. So, if $\sigma$ has support of size $m$, then $A_{\sigma}\cong\Alt(m)$ and $A_{\sigma^\perp}\cong\Alt(n-m)$.

\begin{claim}\label{claim.Aalphaperp}
If $\alpha \in A$ is a bitransposition, then $A_{\alpha^\perp}$ acts faithfully on $C_\alpha$. 
\end{claim}
\begin{proofclaim}
Suppose $A_{\alpha^\perp}$ is not faithful on $C_\alpha$. Let $\{1,\alpha_1,\alpha_2,\alpha_3\}$ be the Klein group on the support of $\alpha$. Since $A_\alpha$ permutes the $C_{\alpha_i}$ transitively while centralizing $A_{\alpha^\perp}$, $A_{\alpha^\perp}$ is not faithful on \emph{any} $C_{\alpha_i}$. 
Since the normal subgroups of $A_{\alpha^\perp}$ are linearly ordered, some nontrivial subgroup of $A_{\alpha^\perp}$ centralizes all $C_{\alpha_i}$, hence all of $G$ by Fact~\ref{fact.4Group}, a contradiction. %(We are also using that $n\ge 7$ implies $A_{\alpha^\perp}$ is nontrivial.)
\end{proofclaim}

\begin{claim}\label{claim.solvable}
If $\alpha \in A$ is a bitransposition, then $C_\alpha$ is solvable but nonabelian.
%In particular, $C_{\alpha,\sigma} \neq 1$ for any involution $\sigma\in A$ that commutes with $\alpha$
\end{claim}
\begin{proofclaim}
Suppose $C_\alpha$ is nonsolvable. We know $A_{\alpha^\perp}\cong \Alt(n-4)$ acts faithfully on $C_\alpha$. By minimality of $n$, we must have $\rk C_\alpha \ge n-4> r-4$, so the corank of $C_\alpha$ in $G$ is at most $3$.   Corollary~\ref{cor.BigCorank} implies $r\le 6$ and then also that $\rk C_\alpha \le 3$, contradicting Fact~\ref{fact.SmallGroups}. Thus $C_\alpha$ is solvable.

We now show  $C_\alpha$ is nonabelian. We may assume $\alpha = (12)(34)$. Let $\beta = (34)(56)$. We claim that $C_{\alpha,\beta}$ is nontrivial. If $C_{\alpha,\beta} = 1$, then $\beta$ inverts $C_\alpha$. The same is then true (by conjugacy under the action of $A_{\alpha^\perp}$) for all bitranspositions $(34)(ij)$ with $5\le i < j$, so $A_{\alpha^\perp}$ centralizes $C_\alpha$. This contradicts  Claim~\ref{claim.Aalphaperp}, so $\rk C_{\alpha,\beta}\ge  1$.  Now, if $C_\alpha$ is abelian, then $C_{\alpha,\beta} = C_{\alpha,\beta,\alpha\beta}$ is a nontrivial central subgroup of $\langle C_{\alpha}, C_{\beta},C_{\alpha\beta}\rangle = G$, against simplicity.
\end{proofclaim}

The next claim handles the generic case, after which we will have only a handful of low rank cases to address.
\begin{claim}\label{claim.CalphaBound}
Let $\alpha=(ab)(cd) \in A$ be a bitransposition and $\Sigma_{\alpha^\perp}$ the subgroup of $A$ fixing $\{a,b\}$ pointwise and $\{c,d\}$ setwise. Then
\[r-6+\ell\le \rk C_\alpha \le r+\frac{1}{2}-\sqrt{2r+\frac{1}{4}}\]
where $\ell$ is the length of any $\Sigma_{\alpha^\perp}$-invariant subnormal series for $C_\alpha$.
In particular, $6\le r\le 10$,
\end{claim}
\begin{proofclaim}
The upper bound  is  Corollary~\ref{cor.BigCorank}. For the lower bound, we apply Lemma~\ref{lem.solvable} to the action of $\Sigma_{\alpha^\perp} \cong \Sym(n-4)$ on $C_\alpha$ (using Claim~\ref{claim.Aalphaperp}) to find that  $\rk C_\alpha \ge (n-4) +\ell-3 = n-7+\ell> r-7+\ell$. Then, since $C_\alpha$ is nonabelian, we have that $\ell \ge 2$, and now the upper and lower bounds together imply that $r\le 10$.  %(Of course, $r\ge 6$ by Lemma~\ref{lem.KleinOnSmall}.)
\end{proofclaim}

Of the low rank cases that remain, $r=6$ stands out for two reasons that quickly arise when considering the action of $A_{\alpha^\perp}$ (or 
$\Sigma_{\alpha^\perp}$) on $C_\alpha$. First, with notation as above, $r=6$ allows for $\ell = 3$ in Claim~\ref{claim.CalphaBound}; whereas, $r\ge 7$ forces  $\ell = 2$. Second, and perhaps more importantly, the $r=6$ case must cover the possibility that $A_{\alpha^\perp}$ is $\Alt(3)$ and hence without involutions.

\begin{claim}
$r=6$.
\end{claim}
\begin{proofclaim}
Assume $7\le r\le 10$. Set $\alpha = (12)(34)$, and let $\Sigma_{\alpha^\perp}$ and $\ell$ be as in Claim~\ref{claim.CalphaBound}. Notice that $\ell\ge 3$ implies $r\le 6$ by Claim~\ref{claim.CalphaBound}, so in the case we are considering, every $\Sigma_{\alpha^\perp}$-minimal subnormal series for $C_\alpha$ has length $2$. We also find  $\rk C_\alpha = r-4$, and (by Corollary~\ref{cor.BigCorank}) $C_\alpha$ is of maximal possible rank among all proper definable subgroups of $G$.

We first show $C_\alpha$ is nonnilpotent. Since $C_\alpha$ is nonabelian, $\rk C_{\alpha,\alpha'}\ge 1$ for $\alpha' = (13)(24)$; also $C_{\alpha,\alpha'} = C_{K}$ for $K = \langle \alpha,\alpha'\rangle$. Now consider $N := N^\circ_{C_{\alpha}}(C_{K})$, and note that $N\neq C_{\alpha}$ as otherwise, by conjugacy under the action of $A_\alpha$, $C_K$ would be  normal in $\langle C_{\alpha}, C_{\alpha'}, C_{\alpha\alpha'}\rangle = G$. Now assume $C_\alpha$ is nilpotent. Then we also have $N \neq C_K$. The series $1 \triangleleft C_{K}\triangleleft N < C_{\alpha}$ is invariant under the action of $\Sigma_{\alpha^\perp}$ (as $\Sigma_{\alpha^\perp}$ centralizes $\alpha$ and normalizes $K$), and by nilpotence of $C_\alpha$, it can be refined to be subnormal and $\Sigma_{\alpha^\perp}$-invariant (and of length at least $3$), a contradiction. 

Thus, we are in the case of $C_\alpha$ being solvable and nonnilpotent. Set $U_\alpha = F^\circ(C_\alpha)$, the connected component of the \emph{Fitting subgroup} of $C_\alpha$; we refer to \cite[Section~7.2]{BoNe94} for the definition and properties of the Fitting subgroup of a group of finite Morley rank. 
Since $C_\alpha$ is nonnilpotent, $1< U_\alpha < C_\alpha$, and as the series has length $2$, it must be $\Sigma_{\alpha^\perp}$-minimal. Observe that $\Sigma_{\alpha^\perp}$ must act faithfully on $U_\alpha$ or  $C_\alpha/U_\alpha$. If not, then $[\beta,\beta,C_\alpha] = 1$ for some bitransposition $\beta \in A_{\alpha^\perp}$, and this implies that the nontrivial elements of $[\beta,C_\alpha]$ (which is nontrivial by Claim~\ref{claim.Aalphaperp}) have order $2$, a contradiction. Thus, by Fact~\ref{fact.CDW}, one of $U_\alpha$ or  $C_\alpha/U_\alpha$ must have rank at least $n-4-2 \ge r-5$, so as $\rk C_\alpha = r-4$, one of the groups has rank $r-5$ and the other rank $1$. We claim that $\rk U_\alpha \neq 1$. 
%Since $\rk C_\alpha \ge 3$, properties of the Fitting subgroup in solvable groups of finite Morley rank  ensure $\rk U_\alpha > 1$. 
This is fairly standard: \cite[Proposition~7.4]{BoNe94} shows $C_\alpha/U_\alpha$ is isomorphic to a quotient of the image of $C_\alpha$ in $\Aut(U_\alpha)$, but if $\rk U_\alpha = 1$, then Zilber's Field Theorem (see~\cite[Theorem~9.1]{BoNe94}) implies that the rank of this image is at most $\rk U_\alpha$, a contradiction since in this case $\rk C_\alpha/U_\alpha \ge 2$.
 We conclude that $A_{\alpha^\perp}$ acts faithfully on $U_\alpha$, which has rank $r-5$ and is an elementary abelian $p$-group (again by Fact~\ref{fact.CDW}).

Let $\beta = (56)(78)$, and note that $W:= \rk C^\circ_{U_\alpha}(\beta)$ has rank at least one. Indeed, if $W$ has rank $0$, then $\beta$  inverts  $U_\alpha$, so as $A_{\alpha^\perp}$ acts faithfully on $U_\alpha$, this would imply that $\beta\in Z(A_{\alpha^\perp}) = 1$. Also, by $p$-unipotence of $U_\alpha$ and solvability of $C_\beta$, we find that $W = (U_\alpha\cap U_\beta)^\circ$. 
Consider $H:=C_G^\circ(W) \ge \langle U_\alpha, U_\beta\rangle$. If $U_\alpha = U_\beta$, then $N_G(U_\alpha) \ge \langle C_\alpha , C_\beta\rangle = G$, a contradiction. Using Corollary~\ref{cor.BigCorank}, we find that $r-5 = \rk U_\alpha < \rk H \le r-4$. As  $U_\alpha$ and $U_\beta$ have corank 1 in $H$, $\rk W = \rk U_\alpha\cap U_\beta = r-6$, so $H/W$ has rank $2$, hence is solvable. Thus, $H$ is solvable, so $H = \langle U_\alpha, U_\beta\rangle$ is nilpotent. We now have $N_G(U_\alpha) \ge \langle C_\alpha, H\rangle = G$, a contradiction. 
\end{proofclaim}

\begin{claim}
$r \neq 6$, and the proof is complete.
\end{claim}
\begin{proofclaim}
 Assume  $r=6$; by Corollary~\ref{cor.BigCorank}, the proper definable subgroups of $G$ have rank at most $3$. 

Let $\alpha = (12)(34)$ with $\Sigma_{\alpha^\perp}$ defined as before. Set $\beta = (34)(56)$, $\gamma = (567)$, and  $K:=\langle\alpha,\beta\rangle$. 
Notice that  conjugacy of $\alpha$, $\beta$, and $\alpha\beta$ implies that $C_{\alpha,\beta} \neq C_\alpha$ as otherwise $C_\alpha = C_\beta = C_{\alpha\beta}$, contradicting Fact~\ref{fact.4Group}. Also notice that $C_{\alpha,\beta}$ is nontrivial since $C_{\alpha}$ is nonabelian (by Claim~\ref{claim.solvable}) and normalized by the involution $\beta$.
Thus, $\rk C_{\alpha,\beta}\ge  1$, which (by our previous observation or Claim~\ref{claim.CalphaBound}) also implies $\rk C_{\alpha}\ge  2$.

We now claim that either $C_\alpha$ or $N_{\alpha,\beta} := N^\circ_G(C_{\alpha,\beta})$ is nonnilpotent. 
Assume $C_\alpha$ is nilpotent (so  $C_\sigma$ is nilpotent for all $\sigma\in K^*$).
The normalizer condition applied in each $C_\sigma$ for $\sigma\in K^*$ ensures that  $\rk( C_\sigma\cap N_{\alpha,\beta})\ge 2$  and, thus, that $\rk N_{\alpha,\beta}$ is at least, hence equal to, $3$. Towards a contradiction, assume $N_{\alpha,\beta}$ is nilpotent. If $\rk C_\alpha = 3$, then the normalizer condition would force $N_{\alpha,\beta}\cap C_\alpha$ to be normal in $G$.  Thus, $\rk C_\alpha=2$, and now the normalizer condition ensures $C_{\alpha,\beta}$ is normal in $C_\sigma$, for  each $\sigma\in K^*$, hence normal in $G$ by Fact~\ref{fact.4Group}, a contradiction since $C_{\alpha,\beta} \neq 1$.

Now consider any proper, definable, connected, nonnilpotent,  $H \le G$; we are thinking of either $C_\alpha$ or $N_{\alpha,\beta}$. Since $H$ has no involutions, the structure of groups of rank at most $3$~\cite{ReJ75,ChG79,FrO18} implies  $H$ must be solvable of rank $3$. An $H$-minimal series for $U:=F^\circ(H)$ must have a factor that is \emph{not} centralized by $H$, and if such a factor has rank $1$, we can linearize the action of $H$ on the factor to produce a rank $1$ field and thus involutions in $H$. Consequently, $U$ must have rank $2$ and be $H$-minimal. Linearizing the action of $H$ on $U$, we find that $T:=H/C_H(U)$ generates a field $L$ in $\End(U)$ and $U$ is $1$-dimensional over $L$. We now specialize this to the case when $H$ is $C_\alpha$ or $N_{\alpha,\beta}$. 

Suppose $C_\alpha$ is nonnilpotent. With notation as above, $U$ and $T$ are $\Sigma_{\alpha^\perp}$-invariant, and the image of $\Sigma_{\alpha^\perp}$ in $\End(U)$ is a group of field automorphisms of $L$. Thus $\Sigma_{\alpha^\perp}$ acts trivially on $L$, hence $L$-linearly on  $U$. Since $U$ is $1$-dimensional, the image of $\Sigma_{\alpha^\perp}$ in $\End(U)$ lies in $L^\times$  forcing $\gamma$ (which is in $ [\Sigma_\alpha^\perp,\Sigma_\alpha^\perp]$) to act trivially on $U$. And considering the action of $\Sigma_{\alpha^\perp}$ on the rank $1$ group $C_\alpha/U$, we find that $[\gamma,\gamma,C_\alpha]=1$. This implies that $[\gamma,C_\alpha]$ is an elementary abelian $3$-group. There are two cases. If $[\gamma,C_\alpha] = 1$, then $C_\alpha \le  C_\gamma$, which we have seen is a contradiction. If $[\gamma,C_\alpha] \neq 1$, then $C_\alpha$-minimality of $U$ implies that $U$ is an elementary abelian $3$-group. In particular, $U$ is normal in every definable connected rank $3$ subgroup containing it. Recall that $U\le C_\gamma$. If equality holds, then conjugacy of $\alpha$ and $(13)(24)$ under the action of $N_A(\langle\gamma\rangle)$ shows $C_\gamma$ is normal in $\langle C_\alpha, C_{(13)(24)} \rangle = G$, a contradiction. And if $U\lneq C_\gamma$, then $U$ is normal in $\langle C_\alpha, C_\gamma \rangle = G$ (since proper  subgroups of $G$ have rank at most $3 = \rk C_\alpha$), again a contradiction.

Finally, suppose $C_\alpha$ is nilpotent; thus $N_{\alpha,\beta}$ is nonnilpotent. We use the notation $U$ and $L$ from above, now applied to $N_{\alpha,\beta}$. This time we find that the image of $K$ in $\End(U)$ lies in $L^\times$, so some nontrivial $\sigma\in K$ centralizes $U$. If $\rk C_\sigma = 3$, then as it is nilpotent, $C_\sigma$ normalizes $U$ and is not equal to $N_{\alpha,\beta}$, forcing $U$ to be normal in $G$. If $\rk C_\sigma = 2$, then $U = C_\sigma\ge C_{\alpha,\beta}$, so as every element of $K$ acts on $U$ as an element of $L^\times$, it must be that $K$ centralizes $U$ (since every element of $K$  centralizes the nontrivial subgroup $C_{\alpha,\beta}$). Thus $U = C_\alpha = C_\beta = C_{\alpha\beta}$, contradicting Fact~\ref{fact.4Group}.
\end{proofclaim}
\end{proof}\setcounter{claim}{0}

% % % % % % % % % % % % % % % % % % % % % % % % % % % % % % % % % % % % % % % % % % % % % % % % % % % % % % % % %
% Section
% % % % % % % % % % % % % % % % % % % % % % % % % % % % % % % % % % % % % % % % % % % % % % % % % % % % % % % % %
\section{A connection to permutation groups}\label{S.PermGroups}

A permutation group $(G, X)$ of finite Morley rank is said to be \emph{generically $t$-transitive} if $G$ has an orbit $\mathcal{O}$ on $X^t$ such that $\rk(X^t \setminus \mathcal{O}) < \rk X^t$.
For example, the natural action of $\GL_n(\bK)$ on $\bK^n$ (with $\bK$ algebraically closed to ensure finiteness of Morley rank) is generically $n$-transitive with  $\mathcal{O}$ being the set of bases. Similarly, $\AGL_n(\bK)$ acts generically $(n+1)$-transitively on $\bK^n$, and  $\PGL_{n+1}(\bK)$ acts generically $(n+2)$-transitively on $\mathbb{P}^n(\bK)$. In fact, in all three examples, every nontrivial element of the group acts without fixing any tuple of the large orbit $\mathcal{O}$; in this case, we say the action is generically \emph{sharply} $t$-transitive.% (generically sharply $(n+2)$-transitive for $\PGL_{n+1}$). 

Notice how generic sharp $t$-transitivity exists (very naturally) for all $t$; this stands in stark contrast to genuine sharp $t$-transitivity, which never exceeds $t=3$ on an infinite set \cite{TiJ52,HaM54}.
Although the degree of generic transitivity is unbounded in general, Borovik and Cherlin showed that there does exists a bound in terms of the rank of the set being acted upon  \cite[Corollary~2.2]{BoCh08}, and they have proposed the problem of showing that there is, in fact, a very natural bound given by (and characterizing!) the $\PGL_{n+1}$ example. To the best of the authors' knowledge, this problem remains open even in the algebraic case, though it has been verified in characteristic  $0$~\cite[Theorem~6.3]{FrMo22}.

\begin{problem*}[{\cite[Problem~9]{BoCh08}}]\label{prob.pgln}
Show that the only transitive and generically $(n+2)$-transitive permutation group of finite Morley rank acting on a set of Morley rank $n$ is $\PGL_{n+1}(\bK)$ acting naturally on the projective $n$-space $\mathbb{P}^n(\bK)$ with $\bK$ algebraically closed.
\end{problem*}

The connection to representations of the symmetric group is as follows.
Suppose $(G,X)$ is generically $t$-transitive, and choose $(x_1,\ldots,x_t)$ in the generic orbit $\mathcal{O}$ of $G$ on $X^t$. Let $G_{t-1} := G_{x_1,\ldots,x_{t-1}}$ be the pointwise stabilizer of $x_1,\ldots,x_{t-1}$, and similarly define $G_{t} := G_{x_1,\ldots,x_{t}}$. Now consider the subgroup $S$ of $G$ that stabilizes $\{x_1,\ldots,x_{t-1}\}$ setwise and fixes $x_t$. It is easily verified that $S$ acts \emph{faithfully} on  $G_{t-1}$, and $S/G_{t}$ is isomorphic to $\Sym(t-1)$ (see \cite[Lemma~4.8, Lemma~4.27]{AlWi18}). In this way, the study of generically $t$-transitive actions is intertwined with faithful representations (on not necessarily abelian groups) of covers of the symmetric group. Moreover, a key point of the problem is that the extreme case of $t\ge \rk X + 2$ should  lead to the identification of $(G,X)$ as  $\PGL_{n+1}(\bK)$ acting naturally on $\mathbb{P}^n(\bK)$, and in that case, $G_{t-1}$ is precisely the maximal torus and $S$ is the Weyl group, further underscoring importance of the subgroups $G_{t-1}$ and $S$.

An important special case of the above situation is the restriction to generically \emph{sharply} $t$-transitive actions. In that case, we have  that $G_{t} = 1$ and $S\cong \Sym(t-1)$. As such, our \hyperlink{h:theorem}{Theorem},  yields the following modest, but nontrivial, corollary.

\begin{maincorollary}
Suppose $(G,X)$ is a generically sharply $t$-transitive permutation group of finite Morley rank with $\rk X=n$. If $t\ge n+2$, then the pointwise stabilizer of a generic $(t-1)$-tuple is solvable or contains involutions.
\end{maincorollary}
\begin{proof}
Assume  $t\ge n+2$; define $G_{t-1}$ and $S\cong\Sym(t-1)$ as in the previous discussion. Then $S$ acts faithfully on $G_{t-1}$. 

We need some additional (and straightforward) background on permutation groups of finite Morley rank. First, $t\ge 2$ implies that  $X$ is connected (i.e.~ of Morley degree 1) \cite[Lemma~1.8(3)]{BoCh08}. Also, the action of $G$  being generically sharply $t$-transitive implies that  $G_{t-1}$ acts generically sharply $1$-transitively on $X$, and this in turn implies  $G_{t-1}$ is in definable bijection with a generic subset of $X$. Thus $G_{t-1}$ is connected of rank  $n$.

Now suppose $G_{t-1}$ is nonsolvable. As connected groups of rank $2$ are solvable~\cite{ChG79}, $n = \rk G_{t-1} \ge 3$, so $t\ge 5$. %(Using Fact~\ref{fact.SmallGroups}, we could even assume $n\ge5$.) 
Consequently, if $G_{t-1}$ does not contain involutions, then everything is in place to apply the \hyperlink{h:theorem}{Theorem} to the action of $S\cong\Sym(t-1)\ge\Sym(4)$ on $G_{t-1}$, but this would imply that  $n \ge t-1$, against our assumption that $t\ge n+2$. 
\end{proof}
The conclusion of the \hyperlink{h:corollary}{Corollary} may be read as saying that either $G_{t-1}$ is ``close'' to being a torus (as desired) or $G_{t-1}$ is subject to further analysis via the rather well-developed theory of groups of finite Morley rank with involutions. Moreover, in the former case,  existing results  (e.g. Lemma~\ref{lem.solvable} and \cite[Theorem]{CDW21})
can be used to further clarify the situation. However, the remaining work needed to address the above problem appears to be significant and surely quite interesting.

% % % % % % % % % % % % % % % % % % % % % % % % % % % % % % % % % % % % % % % % % % % % % % % % % % % % % % % % %
% % % % % % % % % % % % % % % % % % % % % % % % % % % % % % % % % % % % % % % % % % % % % % % % % % % % % % % % %
% % % % % % % % % % % % % % % % % % % % % % % % % % % % % % % % % % % % % % % % % % % % % % % % % % % % % % % % %
% SECTION Acknowledgements
% % % % % % % % % % % % % % % % % % % % % % % % % % % % % % % % % % % % % % % % % % % % % % % % % % % % % % % % %
% % % % % % % % % % % % % % % % % % % % % % % % % % % % % % % % % % % % % % % % % % % % % % % % % % % % % % % % %
% % % % % % % % % % % % % % % % % % % % % % % % % % % % % % % % % % % % % % % % % % % % % % % % % % % % % % % % %
\section*{Acknowledgements}

An early version of this article was discussed with Adrien Deloro, and the authors would like to thank him for his insightful comments and suggestions from which this paper  benefited significantly. The authors also thank the anonymous referee for a careful reading of the paper and many helpful suggestions.

The two authors began discussing this work in January 2021; however, to start in earnest, it was necessary to wait until the first author's return to Lyon in June 2021 after the Turkish state finally agreed to return his passport to him. The first author thanks his coauthor for his frank efforts to keep him active in mathematics. Similar thanks go to Ayşe Berkman, Alexandre Borovik, Gregory Cherlin and Adrien Deloro. 

Both authors warmly thank those who generously invested their efforts in the release of the first author from his two-year-two-month-long forced stay in his homeland which, as of 26 April 2022, has become a forced stay outside his homeland because of a Turkish appeals court decision that again cancelled his passport. Among these companions are many mathematicians. Clearly, mathematics is better done in a more democratic world. 

\bibliographystyle{alpha}
\bibliography{SymBib}
\end{document}